\documentclass[leqno,12pt]{article}
\usepackage{amssymb,amsmath,graphics}
\title{Local Euler-Maclaurin expansion of Barvinok valuations
and Ehrhart coefficients of a rational polytope}
\author{Velleda  Baldoni, Nicole Berline and Mich{\`e}le Vergne}

\date{December 2006}

\begin{document}

\maketitle
\newtheorem{theorem}{Theorem}
\newtheorem{proposition}[theorem]{Proposition}
\newtheorem{lemma}[theorem]{Lemma}
\newtheorem{definition}[theorem]{Definition}
\newtheorem{corollary}[theorem]{Corollary}
\newtheorem{remark}[theorem]{Remark}
\newtheorem{example}[theorem]{Example}
\newenvironment{proof}{{\bf Proof.}}{\par}
\newcommand{\half}{{\frac{1}{2}}}
\newcommand{\C}{{\mathbb C}}
\newcommand{\R}{{\mathbb R}}
\newcommand{\Z}{{\mathbb Z}}
\newcommand{\N}{{\mathbb N}}
\newcommand{\Q}{{\mathbb Q}}
\newcommand{\CA}{{\cal A}}
\newcommand{\CB}{{\cal B}}
\newcommand{\CC}{{\cal C}}
\newcommand{\CE}{{\cal E}}
\newcommand{\CF}{{\cal F}}
\newcommand{\CH}{{\cal H}}
\newcommand{\CL}{{\cal L}}
\newcommand{\CM}{{\cal M}}
\newcommand{\CP}{{\cal P}}
\newcommand{\CS}{{\cal S}}
\newcommand{\CV}{{\cal V}}
\newcommand{\la}{{\langle}}
\newcommand{\ra}{{\rangle}}
\newcommand{\codim}{\operatorname{codim}}
\newcommand{\proj}{\operatorname{proj}}
\newcommand{\Res}{\operatorname{Res}}
\newcommand{\Sum}{\operatorname{Sum}}
\newcommand{\card}{\operatorname{Card}}
\newcommand{\lin}{\operatorname{lin}}
\newcommand{\apex}{\operatorname{apex}}
\newcommand{\relvol}{\operatorname{rvol}}
\newcommand{\Vertex}{\operatorname{Vertex}}
\newcommand{\cone}{\operatorname{cone}}
\newcommand{\vol}{\operatorname{vol}}
\renewcommand{\a}{{\mathfrak{a}}}
\renewcommand{\b}{{\mathfrak{b}}}
\renewcommand{\c}{{\mathfrak{c}}}
\renewcommand{\d}{{\mathfrak{d}}}
\newcommand{\f}{{\mathfrak{f}}}
\renewcommand{\t}{{\mathfrak{t}}}
\newcommand{\p}{{\mathfrak{p}}}
\newcommand{\lattice}{\Lambda}

\section{Introduction}
Let $\p$ be a rational polytope  in $V=\R^d$ and $h(x)$ a
polynomial function on $V$. A classical problem in Integer
Programming is to compute the sum of values of $h(x) $ over the
set of integral points of $\p$,
$$
S(\p,h)=\sum_ {x\in \p\cap \Z^d}h(x).
$$
When $\p$ is dilated by an integer $n\in \N$, we obtain a function
of $n$ which is quasi-polynomial, the so-called Ehrhart
quasi-polynomial of the pair $(\p,h)$
$$
S(n\p,h)= \sum_{m=0}^{d+ N} E_m(\p,h,n)n ^m
$$
of degree $d+N$ where $N=\deg h$. The coefficients $E_m(\p,h,n)$ are
periodic functions of $n\in \N$,  with period the smallest integer
$q$ such that $q\p$ is a lattice polytope.

Replacing $h(x)$ by an exponential, we are led to study the
analytic function on $V^*$ defined by
$$
S(\p)(\xi) = \sum_ {x\in \p\cap \Z^d} e^{\la \xi,x\ra}.
$$
If $\p$ is any rational polyhedron, this sum still makes sense  as a
meromorphic function defined near $0$ and the map $\p\mapsto
S(\p)(\xi)$ is a valuation.

In \cite{EML}, we proved that the meromorphic  function $S(\p)(\xi)$
has a local Euler-Maclaurin expansion
$$
S(\p)(\xi)=\sum_{\f}\mu(\t(\p,\f))(\xi) \int_\f e^{\la \xi,x\ra}dx.
$$

The sum is taken over the set of faces $\f$ of the polyhedron $\p$.
For each face $\f$, the function $\mu(\t(\p,\f))(\xi)$ is
holomorphic near $0$, and it depends only on the transverse cone
$\t(\p,\f)$ of $\p$ along $\f$. More precisely, once a rational
scalar product is chosen on $V$, we define canonically a map
$\a\mapsto \mu(\a)$ from the set of rational affine cones $\a$ in
quotient spaces of $V$,  with values in the space of functions on $V
$ which are analytic near $0$, then we prove that these functions
satisfy the above formula.  The map $\a\mapsto \mu(\a)$ is invariant
under lattice translations, equivariant with respect to lattice
preserving isometries,  and it is a valuation on the set of affine
cones with a fixed vertex (\cite{EML}, Theorems 17 and 18).

It is easy to see that the Ehrhart quasi-polynomial  can be
computed in terms of the Taylor coefficients of the  functions
$\mu(\t(\p,\f))(\xi)$.  For example, if $\p$ is a lattice polytope
and  $h(x)=1$, we have (\cite{EML}, Corollary 28)
\begin{equation}\label{ouralgo}
{\rm Card} (n\p\cap \Z^d) =\sum_{\f} \mu(\t(\p,\f))(0)\vol(\f)
n^{\dim \f}.
\end{equation}

 Using the valuation property of $\mu(\a)$
and
 Barvinok's decomposition of a cone into unimodular cones, we thus
obtained in \cite{EML} an algorithm for computing the Ehrhart
quasi-polynomial. It  has polynomial length with respect to the
input $(\p,h)$, when the dimension $d$ and the degree $N$ are
fixed.

The valuation $S(\p,h)$ has a natural generalization used by
Barvinok in \cite{barvinokEhrhart}, the \emph{mixed} valuation
$S^L(\p,h)$, where $L\subseteq V$ is a rational vector subspace.
 Denote the
projected lattice on $V/L$ by $\lattice_L$.  For a polytope
$\p\subset V$ and a polynomial $h(x)$
$$
S^L(\p,h)= \sum_{y\in \lattice_{V/L}} \int_{\p\cap (y+L)} h(x)dx.
$$
In other words, the polytope $\p$  is sliced along lattice affine
subspaces parallel to $L$ and the integrals of $h$ over the slices
are added up. For $L=V$, there is only one term and $S^V(\p,h)$ is
just the integral of $h(x)$ over $\p$, while, for $L=\{0\}$, we
recover $S(\p,h)$, the sum of values of $h(x) $ over the set of
integral points of $\p$.

In the case $h(x)=1$, we write  $S(\p)$ and $S^L(\p)$ in place of
$S(\p,1)$ and $S^L(\p,1)$.

Using these mixed valuations,  Barvinok gave an algorithm which
computes the $r+1$ highest degree Ehrhart coefficients  of $S(n\p) =
{\rm Card} (n\p\cap \Z^d)$, when $\p$ is a simplex in $\R^d$.
Barvinok's algorithm has polynomial length when $d$ is an input,
provided  $r$ is fixed. The method consists in reducing the problem
to summations over lattice points in dimension $\leq r$.

Barvinok considers particular linear combinations
$$
\sum_{L\in\CL}\rho(L)S^L(\p),
$$
where $\CL$ is a finite set of rational vector subspaces of $V$
which is closed under sum, and the coefficients $\rho(L)$ are
integers which satisfy the following relation between characteristic
functions:
$$
\chi(\cup_{L\in\CL}L^\perp)=\sum _{L\in \CL}\rho(L)\chi(L^\perp),
$$
where $L^\perp\subseteq V^* $ is the orthogonal of $L$. We call a
function $\CL\to \Z$ with this property a patchwork  function on
$\CL$.

When $\p$ is dilated by an integer $n$,  $S^L(n\p)$ is again given
by a quasi-polynomial in $n$, as is a linear combination
$$
\sum_{L\in\CL}\rho(L)S^L(n\p)=\sum _{m=0}^d \nu_m(\p,n)n^m.
 $$

The main theoretical result of \cite{barvinokEhrhart},  Theorem
(1.3), is the following:  if  $\CL$ is  a family of subspaces
which is closed under sum and contains the vector subspace
$\lin(\f)$ parallel to $\f$, for every  face $\f$ of codimension
$\leq r$ of $\p$, and if $\rho$ is a patchwork function on $\CL$,
then the $r+1$ highest degree coefficients $\nu_m(\p,n)$, for $m=
d, \dots, d-r$,  are equal to the corresponding Ehrhart
coefficients of $S(n\p) = {\rm Card} (n\p\cap \Z^d)$.

In the present article, we introduce the meromorphic functions
which extend $S^L(\p)$. For any polyhedron $\p$,
$$
S^L(\p)(\xi)= \sum_{y\in \lattice_{V/L}} \int_{\p\cap (y+L)} e^{\la
\xi,x\ra}dx
$$
is defined as a meromorphic function near $\xi=0$. We show that
$S^L(\p)(\xi)$  also enjoys a local Euler-Maclaurin expansion
(Theorem  \ref{maintheorem}),
$$
S^L(\p)(\xi) =\sum_{\f} \mu^L(\t(\p,\f))(\xi)\int_\f e^{\la
\xi,x\ra}dx.
$$

Furthermore, for a linear combination of Barvinok type, if $\f$ is
a face of $\p$ such that $\lin(\f)\in \CL$, we prove that the
$\f$-term in the Euler-Maclaurin expansions of
$$
S^{\CL,\rho}(\p)(\xi)= \sum_{L\in\CL}\rho(L)S^L(\p)(\xi)
$$
 and of
the \emph{usual} valuation $S(\p)(\xi)$ are equal (Theorem
\ref{fterm}):
$$
\sum_{L\in\CL}\rho(L)\mu^L(\t(\p,\f))(\xi)= \mu(\t(\p,\f))(\xi).
$$
This is the main result of the present article. From the relation
between Ehrhart quasi-polynomials and Euler-Maclaurin expansions,
it implies Barvinok's Theorem (1.3).

Actually, we derive from Theorem \ref{fterm} another computation
of the $r+1$ highest coefficients of the Ehrhart quasi-polynomial
for a pair $(\p,h)$, based on Brion's decomposition of a polytope
into
 cones, in the line of \cite
{barvinok94} and \cite{latte}.

For each vertex $s$ of $\p$,  let $\c_s$ be the cone of feasible
directions of $\p$ at  $s$. Instead of the full family $\CL$
generated by taking sums of the subspaces $\lin(\f)$, when $\f$
runs over the set of
 faces of codimension $\leq r$ of the polytope $\p$, we
consider, for each vertex $s$ of $\p$, the family $\CL_s$
generated by faces of  $\c_s$ of codimension $\leq r$. The point
in taking a family which depends on $s$ lies in the case where
$\p$ is simplicial. Then $\CL_s$ consists only of the spaces
$\lin(\f)$ where $ \f $ is a face of $\c_s $  of codimension $\leq
r$, as this set is already closed under sum. Moreover the
coefficients $\rho(L)$ are just signed binomial coefficients
(Lemma \ref{rho}), and the computation of  $S^L(\c_s)$ is easier
when $L$ is parallel to a  face of $\c_s$ (Example \ref{face}).

Let us describe our method in the simpler case of a  lattice
polytope $\p$ and polynomial $h(x)=1$. By Brion's theorem, we have
$$
S(\p)(\xi)= \sum_s e^{\langle\xi,s\rangle} S(\c_s)(\xi).
$$

For each vertex $s$, let $\rho_s: \CL_s\to \Z$ be a patchwork
function. We define
$$
\CB(\p)(\xi)= \sum_{s}e^{\langle\xi,s\rangle}
S^{\CL_s,\rho_s}(\c_s)(\xi).
$$
For the dilated polytope $n\p$, we have
$$
S(n \p)(\xi) = \sum_{s}e^{n \langle\xi,s\rangle} S(\c_s)(\xi)=
\sum_{m\geq 0}\frac{n^m}{m!}\sum_s \langle \xi,s\rangle^m
\;S(\c_s)(\xi).
$$
Hence, the meromorphic function $\frac{1}{m!}\sum_s \langle
\xi,s\rangle^m \;S(\c_s)(\xi)$ is actually regular at $\xi=0$  and
its value at $\xi=0$ is the $m$th Ehrhart coefficient of $\p$.

 We have similarly
$$
\CB(n\p)(\xi)=  \sum_{m\geq 0}\frac{n^m}{m!}\sum_s \langle
\xi,s\rangle^m \;S^{\CL_s,\rho_s}(\c_s)(\xi).
$$
The meromorphic functions $S^L(\c_s)(\xi)$ and
$S^{\CL_s,\rho_s}(\c_s)(\xi)$ have a special form: they can be
written as the quotient of an analytic function by a product of
$d'\leq d$  linear forms. Such a function $\phi$ has  an expansion
into rational functions $\phi=\sum_{j\geq -d}\phi_{[j]}$ where
$\phi_{[j]}$ is homogeneous of total degree $j$.

Now it follows from our main theorem that, for $m \geq  d-r$, we
have
$$
S(\c_s)_{[-m]}(\xi)= S^{\CL_s,\rho_s}(\c_s)_{[-m]}(\xi),
$$
hence the zero degree terms of $\sum_s \langle \xi,s\rangle^m
\;S(\c_s)(\xi)$ and $\sum_s \langle \xi,s\rangle^m
\;S^{\CL_s,\rho_s}(\c_s)(\xi)$ are equal. Therefore the latter is
also  analytic at $\xi=0$ and its value at $\xi=0$ is the $m$th
Ehrhart coefficient of $\p$. This is the content of Theorem
\ref{betterthanbarvinok}.

Thus,  besides taking care of any polynomial $h(x)$, not only
$h(x)=1$, this method to compute the   $r+1$ highest degree
Ehrhart coefficients for the pair $(\p,h)$ leads to a  simpler
algorithm than the one proposed in \cite{barvinokEhrhart}.   When
$\p$ is  a rational simplex, the contributions of the terms of the
form $S^L(\c_s)(\xi)$ when $L\in \CL_s$  are immediately reduced
to the computation of a function $S(\a)$ with $\a$ a simplicial
cone of dimension smaller or equal to $r$.

There is also another possible implementation of an algorithm to
compute the $r+1$ Ehrhart  highest degree coefficients for the
pair $(\p,h)$ based on the results of \cite{EML}. As seen from
Equation (\ref{ouralgo}), this involves the computation of the
analytic function $\mu(\t(\p,\f))$, also associated  to simplicial
cones in dimension smaller or equal to $r$. We plan to compare the
implementation of  both methods in the near future.

\section{Local Euler-Maclaurin expansion of a mixed valuation $S^L$}

We consider  a rational vector space $V$ of dimension $d$, that is
to say a finite dimensional real vector space with a lattice
denoted by $\lattice_V$ or simply $\lattice$.  We will need to
consider subspaces and quotient spaces of $V$, this is why we
cannot just let $V=\R ^d$ and $\lattice = \Z^d$. By  lattice, we
mean a discrete additive subgroup of $V$ which  generates  $V$ as
a vector space. Hence, a lattice is  generated by a basis of the
vector space $V$. A basis of $V$ which is a $\Z$-basis of
$\lattice_V$ is called an integral basis.
 The elements of $\lattice$
are called integral. An element $x\in V$ is called rational if
$qx\in \lattice $ for some integer $q\neq 0$. The space of
rational points in $V$ is denoted by $V_\Q$. A subspace $L$ of $V$
is called rational if $L\cap \lattice $ is a lattice in $L$.  If
$L$ is a rational subspace, the image of $\lattice$ in $V/L$ is a
lattice in $V/L$, so that $V/L$ is a rational vector space. We
will call the image of $\lattice$ in $V/L$ the projected lattice.
\begin{example}\label{projectedlattice}
Let $V=\R^2$ with standard lattice $\Z^2$. Let $v_1, v_2$ be two
primitive integral independent vectors. Using an integral  basis
with first basis vector  $v_1$,  a straightforward computation
shows that the projected lattice on $\R^2 /\R v_1$  is $\Z
\frac{{\overline v_2}}{\det(v_1,v_2)} $, where ${\overline v_2}$
is the projection of $v_2$ on $\R^2 /\R v_1$.
\end{example}

A rational space $V$, with lattice $\lattice$, has a canonical
Lebesgue measure, for which  $V/\lattice$ has measure $1$. An affine
subspace $L$ of $V$ is called rational if it is  a translate of a
rational subspace by a rational element. It is similarly provided
with a canonical Lebesgue measure. We will  denote this measure  by
$dm_L$.

We will denote elements of $V$ by latin letters $x,y,v,\dots$ and
elements of $V^*$ by greek letters $\xi,\alpha,\dots$.  We denote
the duality bracket by $\langle\xi,x\rangle$.

If $S$ is a subset of $V$, we denote by $<S>$ the affine subspace
generated by $S$.  If  $S$ consists of rational points, then
$<S>$ is rational. Remark that  $<S>$ may  contain no integral
point.
 We denote by $\lin(S)$ the vector subspace of $V$ parallel to
$<S>$.

If $S$ is a subset of $V$, we denote by $S^{\perp}$ the subspace of
$V^*$ orthogonal to $S$:
$$
S^{\perp}= \{\xi\in V^*\; ;\,\,\langle\xi,x\rangle =0 \;\mbox{for
all}\;
  x\in  S\}.
$$
If $L$ is a subspace of $V$, the dual space  $(V/L)^*$ is
canonically identified with the subspace $L^{\perp} \subset V^*$.

A convex rational polyhedron $\p$ in $V$ (we will simply say
 polyhedron) is, by definition, the intersection of a finite number of
half spaces bounded by  rational affine hyperplanes. We say that
 $\p$ is solid (in V) if $<\p>= V$.  A polytope $\p$ is a {\bf compact} polyhedron.

The set of non  negative real numbers is denoted by $\R_+$. A convex
rational cone
 $\c$ in $V$ is  a closed
convex cone $\sum_{i=1}^k\R_+ v_i$ which is generated
 by a finite number of elements  $v_i$ of $V_\Q$.
In this article, we simply say  cone instead of convex rational
cone.

An affine (rational) cone $\a$ is, by definition, the translate of a
cone in $V$ by a rational element $s\in
  V_\Q$. This cone is uniquely defined by $\a$.

 A  cone $\c$ is called simplicial if it is  generated by
independent elements of $V_\Q$. A simplicial cone $\c$ is called
unimodular if it is generated by independent integral vectors
$v_1,\dots, v_k$ such that $\{v_1,\dots, v_k\}$ can be completed
in an integral basis of $V$. An affine cone $\a$ is called
simplicial (resp. simplicial unimodular) if it is the translate of
a simplicial (resp. simplicial unimodular) cone.

An affine cone $\a$ is called pointed if it does not contain any
straight line.

The set of faces of an affine cone $\a$ is denoted by $\CF(\a)$. If
$\a$ is pointed, then the vertex of $\a$ is the unique face of
dimension $0$, while $\a$ is the unique face of maximal dimension
$\dim \a$.

 Let us recall the definition of the \emph{transverse cone} $\t(\p,\f)$
of a polyhedron $\p$ along one of its faces $\f$. Let $x$ be a
point in the relative interior of $\f$. The cone of feasible
directions of $\p$  at $x$ is the set $\c(\p,\f):=\{v\in V\;;
x+\epsilon v\in \p \;\mbox{for}\; \epsilon>0\; \mbox{small enough}
\}$. It does not depend on the choice of  $x$. We denote the
projection $V\to V/\lin(\f)$ by $\pi_\f$. Then $\t(\p,\f)$ is the
image  $ \pi_\f(\f+\c(\p,\f))$ of the affine
 cone $\f+\c(\p,\f)$ in  $V/\lin(\f)$. It is a solid
pointed \textbf{affine} cone in the quotient space $V/\lin(\f)$
with vertex $\pi_\f(<\f>)$. In particular, if $v$ is a vertex of
$\p$, the transverse cone $\t(\p,v)$ coincides with the supporting
cone $v+\c(\p,v)\subset V
 $.

If $\a$ is an affine cone and $\f$ is a face of $\a$,  then
$\c(\a,\f)= \a + \lin(\f)$ and the transverse cone $\t(\a,\f)$ of
$\a$ along  $\f$ is just  the projection $\pi_\f(\a)$  of $\a$ on
$V/\lin(\f)$.

\begin{definition}
Denote by  $\CH(V^*)$  the ring of analytic functions around $0\in
V^*$. Denote by  $\CM(V^*)$ the ring of meromorphic functions
defined around $0\in V^*$  and by $\CM_{\ell}(V^*)\subset
\CM(V^*)$ the subring consisting of those meromorphic functions
$\phi(\xi)$ such that there exists a product of linear forms
$D(\xi)$ with
$$
D(\xi)\phi(\xi)\in \CH(V^*).
$$
A function $\phi(\xi)\in \CM_{\ell}(V^*)$ has a unique expansion
into homogeneous rational functions
$$
\phi(\xi)= \sum_{m>>-\infty}\phi_{[m]}(\xi)
$$
where $m$ is the total degree.
\end{definition}

If $P$ is a homogeneous polynomial on $V^*$ of degree $p$, and $D$
a product of $r$ linear forms, then $\frac{P}{D}$ is an element in
$\CM_{\ell}(V^*)$  homogeneous of degree $m=p-r$.

Let us recall the definition of the function $I(\p)\in
\CM_{\ell}(V^*)$ associated to a polyhedron $\p$, (see for
instance the survey \cite{barpom}).
\begin{proposition}\label{valuationI}
There exists a map $I$ which  to every polyhedron $\p\subset V$
associates a meromorphic function with rational coefficients
$I(\p)\in \CM_{\ell}(V^*)$, so that the following properties hold:

(a) If $\p$ contains a straight line, then $I(\p)$=0.

(b) If $\xi\in V^*$ is such that $|e^{\la \xi,x\ra}|$ is integrable
over $\p$,  then
$$
I(\p)(\xi)= \int_\p e^{\la \xi,x\ra} dm_{<\p>}(x).
$$

(c) For every point $s\in V_\Q$, we have
$$
I(s+\p)(\xi) = e^{\la \xi,s\ra}I(\p)(\xi).
$$

(d) The map $I$ is a {\em simple valuation}: if the characteristic
functions $\chi(\p_i)$ of a family of  polyhedra $\p_i$ satisfy a
linear relation $\sum_i r_i \chi(\p_i)=0$, then the functions
$I(\p_i)$ satisfy the  relation
$$
\sum_{\{i,<\p_i>=V\}}r_i I(\p_i)=0.
$$
\end{proposition}
\medskip

In the following proposition, we define the \emph{mixed valuation}
$\p\mapsto S^L(\p)$ associated to a rational vector subspace
$L\subseteqq V$. To any polyhedron $\p$, we associate a
\emph{meromorphic function} $S^L(\p)(\xi)\in \CM(V^*)$. If $\p$ is
compact, this function is actually regular at $0$, and its value
for $\xi=0$ is the valuation $E_{L^\perp}(\p)$ considered by
Barvinok \cite{barvinokEhrhart}.

We denote by $\lattice_{V/L}$ the projection on $V/L$ of the lattice $\lattice$.

\begin{proposition}\label{valuationSL}
Let $L\subseteqq V$ be a rational subspace. There exists a map $S^L$
which to every rational polyhedron $\p\subset V$ associates a
meromorphic function with rational coefficients $S^L(\p)\in
\CM(V^*)$ so that the following properties hold:

(a) If $\p$ contains a line, then $S^L(\p)$=0.

(b)
\begin{equation}\label{SL}
S^L(\p)(\xi)= \sum_{y\in \lattice_{V/L}} \int_{\p\cap (y+L)}
e^{\la \xi,x\ra}dm_L(x),
\end{equation}
for every  $\xi\in V^*$  such that the above sum converges.

(c) For every  point $s\in \lattice$, we have
$$
S^L(s+\p)(\xi) = e^{\la \xi,s\ra}S^L(\p)(\xi).
$$

(d) The map $S^L$ is a valuation: if the characteristic functions
$\chi(\p_i)$ of a family of  polyhedra $\p_i$ satisfy a linear
relation $\sum_i r_i \chi(\p_i)=0$, then the functions $S^L(\p_i)$
satisfy the same relation
$$
\sum_i r_i S^L(\p_i)=0.
$$
\end{proposition}
For $L=\{0\}$, we recover the valuation $S$ given by
$$
S(\p)(\xi)= \sum_{x\in \p\cap \lattice}e^{\la \xi,s\ra},
$$
provided this sum is convergent.

 For $L=V$, we have $S^V(\p)=I(\p)$,
if $\p$ is solid, and $S^V(\p)=0$ otherwise.

The proof is  entirely analogous to the case  $L= \{0\}$, see
Theorem 3.1 in \cite{barpom}, and we omit it.
\begin{remark} The function
$S^L(\p) $ is actually an element of $\CM_{\ell}(V^*)$, but we do
not prove it at this point. Let $\a$ be an affine cone and
$\{v_i\}$ the generators of its edges. It will follow from the
Euler-Maclaurin expansion of $S^L(\a) $ (Theorem
\ref{maintheorem})  that  $\prod_i \langle \xi,v_i\rangle
S^L(\a)(\xi)$ is analytic near zero for any $L$. It would be
interesting to prove it a priori. By Brion's theorem and the
valuation property, it follows in particular that $S^L(\p) \in
\CM_{\ell}(V^*)$.
\end{remark}
\begin{example}\label{face}
Let $\a$ be a simplicial  affine cone in the space $V$, and assume
that $L=\lin(\f_1)$ for some face $\f_1$ of $\a$. In this case,
$S^L(\a)(\xi)$ decomposes as product of an integral and a discrete
sum. For simplicity, assume that $\a$ is solid. Let $\f_2 $  be
the  face of $\a$ such  $V= \lin(\f_1)\oplus \lin(\f_2) $. We
write $x=x_1 +x_2 $ and $\xi= \xi_1+ \xi_2$ for the corresponding
decompositions of $x\in V$ and $\xi\in V^*$. Then $\a=\a_1+\a_2$
where $\a_i$ is a simplicial affine cone in $\lin(\f_i)$.  Let us
denote by $\lattice_2$ the projected lattice in $V/\lin(\f_1)\sim
\lin(\f_2)$. From (\ref{SL}), we obtain immediately
\begin{equation}\label{sumface}
S^L(\a)(\xi_1+ \xi_2)= I(\a_1)(\xi_1) \sum_{x_2\in \a_2\cap
\lattice_{2}}e^{\langle \xi_2,x_2\rangle}.
\end{equation}
Notice that the lattice $\lattice_2$ is usually bigger than
$\Lambda\cap \lin(\f_2)$.
\end{example}
\begin{example}\label{dim2}
Let $V=\R^2$ with the standard lattice. We compute $S^L(\a)$ when
$\a$ is a cone and $L$ is a line. Let $\a= \R_+ v_1+ \R_+ v_2$,
where  $v_1, v_2$ are two linearly independent primitive integral
vectors.

(a)  Assume that $L$ is the line supporting an edge of $\a$, say $L=
\R v_1$. We identify $V/L$ to $\R v_2$. The projected lattice is
$\lattice_2= \Z \frac{v_2}{\det(v_1,v_2)} $, (Example
\ref{projectedlattice}), hence, by (\ref{sumface}) in Example
\ref{face}, we have
\begin{equation}\label{dim2a}
S^L(\a)(\xi)=-\frac{1}{\langle \xi,v_1\rangle}\; \frac{1}{
1-e^{\frac{\langle \xi,v_2\rangle}{|\det(v_1,v_2)|}}}.
\end{equation}

(b) Assume now that $L$ is transverse to both edges of $\a$.
Assume that $ \det(v_1,v_2)>0$. Let $L=\R u$ where $u$ is a
primitive integral vector chosen so that $\det(u,v_2)>0$. Let
$\a_i= \R_+u + \R_+ v_i$ for $i=1,2$. We decompose the
characteristic function of the cone $\a$ as $\chi(\a)= \chi(\a_2)
+ \chi(\a_1)- \chi(\R_+ u)$ or $\chi(\a)= \chi(\a_2)-\chi(\a_1)+
\chi(\R_+ v_1)$, depending on whether $u$ belongs to $\a$  or not.
Using the valuation property, case {\it (a)} and the relation
$$
\frac{1}{1-e^x}+\frac{1}{1-e^{-x}}=1,
$$
we obtain in both cases
\begin{equation}\label{dim2b}
S^L(\a)(\xi)=-\frac{1}{\langle \xi,u\rangle}\left(
\frac{1}{1-e^{\frac{\langle \xi,v_2\rangle}{\det(u,v_2)}}}-
\frac{1}{1-e^{\frac{\langle \xi,v_1\rangle}{\det(u,v_1)}}}
\right).
\end{equation}

In this example, we see that $\langle \xi,v_1\rangle \langle
\xi,v_2\rangle S^L(\a)(\xi)$ is indeed analytic near $\xi=0$.
\end{example}

In the following theorem and its applications, we will consider
the functions $S^L(\p)$ when the space $V$ is replaced with a
quotient space $W$. We denote by $\CC(W)$ the set of affine cones
in $W$. Thus if $\a \in \CC(W)$, and $L$ a rational subspace of
$W$, the function $S^L(\a)$ is a meromorphic function on $W^*$. We
are going to show that the function $S^L(\a)$ has a local
Euler-Maclaurin expansion, which generalizes the case $L=\{0\}$ of
\cite{EML}.

\begin{theorem}\label{maintheorem}
Let $V$ be a rational space and $Q$ a rational scalar product on
$V^*$. There exists a unique family of maps  $\mu_W^L$, indexed by
pairs $(W,L)$ where $W$ is a rational quotient space   of $V$ and
$L$ is a rational vector subspace of $W$ such that the family
enjoys the following properties:

(a) $\mu_W^L$ maps  $\CC(W)$ to  $\CH(W^*)$, the space of analytic
functions on $W^*$.

(b) If $W=\{0\}$,  then  $\mu_{\{0\}}^{\{0\}}(\{0\})= 1 $.

(c) For $\dim W >0$ and $L=W$, then $\mu_W^W(\a)=0$.

(d) If the affine cone $\a\in \CC(W)$ contains a straight line,
then $\mu_W^L(\a)=0$.

(e) For any affine cone $\a$ in $W$, the following formula holds
$$
S^L(\a)=
\sum_{\f\in\CF(\a)}\mu_{W/\lin(\f)}^{L+\lin(\f)/\lin(\f)}(\t(\a,\f))I(\f)
$$
where the sum is over all faces of the cone $\a$.

\bigskip

In this last formula, the function
$\mu_{W/\lin(\f)}^{L+\lin(\f)/\lin(\f)}(\t(\a,\f))$ is considered
as a function on $W^* $ itself by means of the orthogonal
projection $W^* \to (W/\lin(\f))^*= (\lin(\f))^\perp $ with
respect to the scalar product on $W^*\subset V^*$.
\end{theorem}

\begin{proof}
The proof is entirely similar to the case  $L=\{0\}$ studied in
\cite{EML}.  Note that $\mu_W^{\{0\}}$ coincides with the map
denoted by $\mu_W$ in \cite{EML}. The only new item is (c). It
follows immediately from the relation $S^W(\a)= I(\a)$.
\end{proof}

\begin{remark}
Let $\a$ be a solid cone in $W$, and let $\f$ be a face of $\a$
such that $\dim\f< \dim W$. If $L$ is transverse to the face $\f$,
that is, if $L+\lin(\f)=W$, then
$\mu_{W/\lin(\f)}^{L+\lin(\f)/\lin(\f)}(\t(\a,\f))=0$. This
follows from (c).
\end{remark}

From now on we omit the subscript $W$, thus we write $\mu^{L}$ in
place of $\mu_W^L$. The next theorem and its proof are also
entirely similar to the case $L=\{0\}$ in \cite{EML}.
\begin{theorem}\label{maintheoremplusmuL}
The analytic functions  defined in Theorem \ref{maintheorem} have
the following properties:

(a) For any $x\in \lattice$, one has $\mu^L(x +\a)= \mu^L(\a)$.

(b) The map $(\a,L)\mapsto \mu^L(\a)$  is equivariant with respect
to lattice-preserving  isometries. In other words, let $g$ be an
isometry of $W$ which preserves the lattice $\lattice_W$. Then
$\mu^{g(L)}(g(\a))( ^{t}g^{-1}\xi)=\mu^L(\a)(\xi)$.

(c) If $W$ is an orthogonal sum $W=W_1\oplus W_2$ of rational
spaces, $L_i\subseteqq W_i$ and $\a_i$ is an affine cone in $W_i$
for $i=1,2$, then
$$
\mu^{L_1\oplus L_2}(\a_1 + \a_2)= \mu^{L_1}(\a_1)\mu^{L_2}(\a_2).
$$

(d)  For a fixed $s\in W_\Q$, the map $\c\to \mu^L(s+\c)$ is a
valuation on the set of cones in $W$.

(e)  Let $\p\subset W$ be a rational polyhedron, then
\begin{equation}
S^L(\p)=
\sum_{\f\in\CF(\p)}\mu^{L+\lin(\f)/\lin(\f)}(\t(\p,\f))I(\f).
\end{equation}
\end{theorem}

\begin{example} Let us compute the function $\mu^L$ for the various
transverse cones of Example \ref{dim2}. We define a function $B(u)$
on $\C$, holomorphic near $0$, by
$$
B(u)=\frac{1}{1-e^u}+\frac{1}{u}.
$$
We have
$$
I(\a)(\xi)= \frac{|\det(v_1,v_2)|}{\langle \xi,v_1\rangle \langle
\xi,v_2\rangle}.
$$
Consider case (b) where $L$ is transverse to both edges
 $\f_i=\R_+ v_i$ of $\a$.

Using the equation $\det(v_1,v_2)u=\det(u,v_2)v_1-\det(u,v_1)v_2$,
we have
$$
I(\a)(\xi)= \frac{1}{\langle \xi,u\rangle}
\left(\frac{\det(u,v_2)}{\langle \xi,v_2\rangle} -
\frac{\det(u,v_1)}{\langle \xi,v_1\rangle}\right).
$$

Thus we can rewrite (\ref{dim2b}) as
$$
S^L(\a)=\mu^L(\a) +
\sum_{i=1,2}\mu^{L+\lin(\f_i)/\lin(\f_i)}((\t(\a,\f_i))I(\f_i) +
I(\a),
$$
with
\begin{eqnarray}\label{muLdim2transverse}
&&\mu^L(\a)(\xi)= \frac{1}{\langle \xi,u\rangle}
\left[B\left(\frac{\langle
\xi,v_1\rangle}{\det(u,v_1)}\right)-B\left(\frac{\langle
\xi,v_2\rangle}{\det(u,v_2)}\right)\right],
\\
\nonumber && \mu^{L+\lin(\f_i)/\lin(\f_i)}((\t(\a,\f_i))=0 \;\;
{\rm for}\; i=1,2,
\end{eqnarray}
Observe that (\ref{muLdim2transverse}) is indeed regular at $\xi=0$.

In case {\it (a)} where $L= \R v_1$, we have
$L+\lin(\f_1)/\lin(\f_1)=\{0\}$. Let us assume that
$\det(v_1,v_2)>0$. Then we have, by \cite{EML},
$$
\mu^{\{0\}}(\t(\a,\f_1))(\xi)= B\left(\frac{- C_1\langle
\xi,v_1\rangle+\langle \xi,v_2\rangle}{\det(v_1,v_2)}\right)
$$
with $C_1= \frac{Q(v_1,v_2)}{Q(v_1,v_1)}$.

As $ I(\f_1)(\xi)=- \frac{1}{\langle \xi,v_1\rangle} $,  we can
rewrite (\ref{dim2a}) as
$$
S^L(\a)=\mu^L(\a)+ \mu^{\{0\}}(\t(\a,\f_1))I(\f_1)+\;
\mu^{L+\lin(\f_2)/\lin(\f_2)}((\t(\a,\f_2))I(\f_2)+ I(\a),
$$
with
$$
\mu^L(\a)(\xi)=\frac{1}{\langle \xi,v_1\rangle}
\left[B\left(\frac{- C_1\langle \xi,v_1\rangle+\langle
\xi,v_2\rangle}{\det(v_1,v_2)}\right) -B\left(\frac{\langle
\xi,v_2\rangle}{\det(v_1,v_2)}\right)\right],
$$
which is indeed regular at $\xi=0$, and, again,
$$
\mu^{L+\lin(\f_2)/\lin(\f_2)}((\t(\a,\f_2))=0.
$$
\end{example}

\section{Barvinok valuations}
Let $\CL$ be a finite family of rational vector subspaces of $V$,
and let $\rho(L), L\in \CL$, be a set of rational coefficients.
The linear combination $\sum_{L\in \CL} \rho(L) S^L(\p)$  is again
a valuation on the set of polyhedra $\p\subset V$, with values in
$\CM_{\ell}(V^*)$. By taking linear combinations, we obtain a
local Euler-Maclaurin expansion for the function $\sum_{L\in \CL}
\rho(L) S^L(\p)(\xi)$.
\begin{definition}Let $\p\subset V$ be a polyhedron.

(a) We denote
$$ S^{(\CL,\rho)}(\p)=\sum_{L\in \CL} \rho(L) S^L(\p).
$$
(b) We define the $\f$-term in the local Euler-Maclaurin expansion
of $S^{(\CL,\rho)}(\p)$ to be
$$
\mu^{(\CL,\rho)}(\t(\p,\f))=\sum _{L\in\CL}\rho_{\CL}(L)
\mu^{L+\lin(\f)/\lin(\f)}(\t(\p,\f)).
$$
\end{definition}
Thus we have
$$
S^{(\CL,\rho)}(\p)=
\sum_{\f\in\CF(\p)}\mu^{(\CL,\rho)}(\t(\p,\f))I(\f).
$$
We are going to compute the $\f$-term in the case of the following
particular linear combinations introduced by Barvinok
\cite{barvinokEhrhart}.

\begin{definition}
 The valuation  $S^{(\CL,\rho)}$ is called a Barvinok valuation if

(a) the family of subspaces $\CL$ is stable under sum,

(b) $\rho$ is an integer valued function on the set $\CL$ such that
the characteristic function of the union of the subspaces $L^\perp
\subseteqq V^*$ can be written as the linear combination
\begin{equation}\label{rhoL}
 \chi(\cup_{L\in\CL}L^\perp)=\sum _{L\in
\CL}\rho(L)\chi(L^\perp).
\end{equation}
\end{definition}
\begin{definition}
We call a   function $\CL \to \Z$ which satisfies (\ref{rhoL}) a
patchwork function on $\CL$.
\end{definition}

As the set of orthogonal subspaces $L^\perp \subseteqq V^*$ is
stable under intersection, a particular function $\rho_{\CL}$ with
this property can be computed in terms of the Moebius function of
the partially ordered set $\CL$ (\cite{Stanley}, vol I, section
3.7), as explained in \cite{barvinokEhrhart}.

Let us compute a patchwork function $\rho$ in the following case.
$V=\R^d$ with standard basis $e_i, i= 1,\dots,d$, and $\CL_{d,q}$
is the set of subspaces
 $L_I= \oplus _{i\in I}\R e_i$ with cardinal $|I|\geq q$. The function $\rho_{d,q}$
 defined below is actually the one associated
to the Moebius function, but we will not need this fact.

We denote the binomial coefficient $\frac{m!}{k! (m-k)!}$  by
$C^{k}_{m}$.

\begin{lemma}\label{rho}
The function $\rho_{d,q}$ on $\CL_{d,q}$ defined by
$$
\rho_{d,q}(L_I)= (-1)^{n-q}C^{q-1} _{n-1} \;\;\; \mbox{ if}\;\;
|I|=n,
$$
satisfies Equation (\ref{rhoL}).
\end{lemma}

\begin{proof}
If $e^i$ is the dual basis, the orthogonal space $L_I^{\perp}$ is
equal to $\sum_{i\notin I}\R e^i$.
 Let
$\xi=\sum_{i=1}^d \xi_i e^i \in\cup_{L\in\CL_q}L^\perp$.  Let
$I_0$ be the set of indices $i\in [1,\dots, d]$ such that $\xi_i=
0$. Then $|I_0|\geq q$, and  $\xi\in L_I^{\perp}$ if and only
$I\subseteq I_0$. Let $|I_0|=N$. The value at $\xi$ of the
right-hand side of (\ref{rhoL}) is equal to
$$
E(N,q)=\sum_{I\subseteq
I_0}\rho_{d,q}(L_I)=\sum_{n=q}^{N}(-1)^{n-q} C_N^n
 C_{n-1}^{q-1}.
$$
We want to prove that $E(N,q)= 1$. Writing $n=q+i$, we have
$$
E(N,q)
=\sum_{i=0}^{N-q}(-1)^{i}\frac{N!}{(q+i)(N-q-i)!i!(q-1)!}.$$ Let
us compute $(q-1)! \Big(E(N+1,q)-E(N,q)\Big)$. This is equal to
$$
(-1)^{N+1-q}\frac{(N+1)!}{(N+1)(N+1-q)!}
+\sum_{i=0}^{N-q}(-1)^{i}\frac{1}{(q+i)i!}\Big(\frac{(N+1)!}{(N+1-q-i)!}-
\frac{N!}{(N-q-i)!}\Big)$$
$$=N! \sum_{i=0}^{N+1-q}(-1)^{i}\frac{N!}{i!(N+1-q-i)!}=0.$$
We obtain $E(N,q)=E(q,q)=1$. $\Box$
\end{proof}

\bigskip

In the case of a Barvinok valuation, it turns out that the $\f$-term
in the Euler-Maclaurin expansion of
 $S^{(\CL,\rho)}(\p)$ coincides with that of $S(\p)$, if the vector
subspace $\lin(\f)$ belongs to the set $\CL$. This is the crucial
result of the present article. It is an easy consequence of the
following combinatorial lemma.
\begin{lemma}\label{combinatorial}
Let $\cal{L}$ be a finite family of vector subspaces of  $V$, stable
under sum and let $\rho$ be a patchwork function on $\cal{L}$.

(a) Let $L_0\in \CL$. Then
$$ \sum _{\{L\in \CL, \; L\subseteq
L_0\}}\rho(L)=1.
$$

(b) Let $L_0\subsetneqq L_1$ be two  subspaces in the family $\CL$.
Then
$$
\sum _{\{L\in \CL, \; L+L_0=L_1\}} \rho(L) =0.
$$
\end{lemma}
\begin{proof}
There exists a $\xi_0\in L_0^\perp$ such that, for  $L\in \CL$,
$\xi_0 \in L^\perp$ if and only if $L\subseteq L_0$. We obtain
{\it (a)} by evaluating both sides of  (\ref{rhoL}) at this
particular element $\xi_0$. Next, we deduce  {\it (b)} from {\it
(a)}, by induction on $\dim L_1- \dim L_0$.

For two subspaces $M\subseteq M'$ in the family $\CL$,  let us
denote
$$f(M,M')=\sum _{\{L\in
\CL, \; L+M=M'\}} \rho(L).
$$
If $M=M'$,  we have $f(M,M)= \sum _{\{L\in \CL, L\subseteq
M\}}\rho(L)=1$ by {\it (a)}.

We apply this with $M=L_1$. Thus
$$
\sum _{\{L\in \CL, L\subseteq L_1\}}\rho(L)= 1.
$$

In this sum, we group the $L\in\CL$ such that $L+L_0$ is equal to
a given $M_1\in \CL$ together.

First we consider the case when $\dim L_1- \dim L_0= 1$. Then $M_1$
is either $L_0$ or $L_1$, hence we obtain
$$
 f(L_0,L_0) + f(L_0,L_1)= 1.
$$
Since $f(L_0,L_0)=1$ by {\it (a)}, we obtain $f(L_0,L_1)= 0$ as
required.

Next we consider the  case when $\dim L_1- \dim L_0 > 1$. We obtain
$$
\sum _{\{M_1\in \CL, \; L_0\subseteq M_1 \subseteq
L_1\}}f(L_0,M_1)= 1.
$$
For $M_1=L_0$,  we have $f(L_0,L_0)=1$ by {\it (a)}. For
$M_1\subsetneqq L_1$,  we have $f(L_0,M_1)= 0$ by induction
hypothesis. Hence there remains only the term $f(L_0,L_1)$ which
must be equal to $0$. $\Box$
\end{proof}
\medskip

We study now the Euler-Maclaurin expansion of a Barvinok valuation.

\begin{theorem}\label{fterm} Let $\p\subset V$ be a rational polyhedron  and let $\f$ be a face of $\p$.
Let $\CL$ be a finite family of rational vector subspaces of  $V$,
stable under sum. Let $ \rho$ be a patchwork function on $\CL$,
 and let $S^{(\CL,\rho)}=\sum_{L\in\CL}\rho(L)S^L$.

Assume that $\lin(\f)$ belongs to $\CL$. Then
$$
\mu^{(\CL,\rho)}(\t(\p,\f))= \mu^{\{0\}}(\t(\p,\f)).
$$
In other words,  the $\f$-term in the local Euler-Maclaurin
expansion of $S^{(\CL,\rho)}(\p)$ coincides with that of $S(\p)$.
\end{theorem}
\begin{proof}
 In the sum $\sum _{L\in\CL}\rho(L) \mu^{L+\lin(\f)/\lin(\f)}(\t(\p,\f))$,
 we group the terms for which $L+\lin(\f)$ is equal to a given $L_1$ together.
 By Lemma \ref{combinatorial} we obtain $\mu^{\{0\}}(\t(\p,\f))$ for $L_1= \lin(\f)$ and $0$ otherwise. $\Box$
\end{proof}

\begin{corollary}\label{lowestdegree}
Let $\p\subset V$ be a rational polyhedron. Let $0\leq k\leq d$.
Let $\CL$ be a finite family of rational vector subspaces of  $V$,
stable under sum, such that $\lin(\f)\in \CL $ for every
$k$-dimensional face $\f$ of $\p$ and let $ \rho$ be a patchwork
function on $\CL$.

$\bullet$ Let $0< k\leq d$. Then the meromorphic function
$$
S(\p)(\xi)- S^{(\CL,\rho)}(\p)(\xi)
$$
 has lowest degree $\geq -k+1$.

$\bullet$ Let $k=0$. Then
$$ S(\p)(\xi)= S^{(\CL,\rho)}(\p)(\xi).
$$
\end{corollary}

\begin{proof}
By Theorem \ref{fterm}, the local Euler-Maclaurin expansion of the
difference involves only faces of dimension $<k$.
$$
S(\p)(\xi)- S^{(\CL,\rho)}(\p)(\xi)= \sum_{\{\f\in\CF(\p),\; \dim\f
<k\}}\left(\mu^{\{0\}}(\t(\p,\f))(\xi)-
\mu^{(\CL,\rho)}(\t(\p,\f))(\xi)\right)I(\f)(\xi).
$$
For a face of dimension $j$, the function $I(\f)(\xi)$ is
homogeneous of degree $-j$. Multiplied by the  \emph{holomorphic}
function $\mu^{\{0\}}(\t(\p,\f))(\xi)- \mu^{\CL}(\t(\p,\f))(\xi)$,
its lowest degree can only increase.  $\Box$\end{proof}

\bigskip

Remark that the statement of Corollary \ref{lowestdegree} above
does not involve any scalar product. In the next section, we will
show that this corollary implies our main Theorem (Theorem
\ref{betterthanbarvinok}).

\section{Ehrhart quasi-polynomial}
 Let $\p$ be a rational polytope and let $h(x)$ be a polynomial function on $V$.
Let
$$
S(\p,h)=\sum_ {x\in \p\cap \Z^d}h(x).
$$
When $\p$ is dilated by a non negative  integer $n$, we obtain the
Ehrhart quasi-polynomial of the pair $(\p,h)$
$$
S(n\p,h)= \sum_{m=0}^{d+ N} E_m(\p,h,n)n ^m,
$$
where $N=\deg h$. The coefficients $E_m(\p,h,n)$ are periodic
functions of $n\in \N$,  with period the smallest integer $q$ such
that $q\p$ is a lattice polytope.

If an integer $r\leq d$ is fixed, and $h=1$, Barvinok
\cite{barvinokEhrhart} proved that the $r+1$ highest Ehrhart
coefficients $E_d(\p,1,n), \ldots, E_{d-r}(\p,1,n)$ of $S(n\p,1)$
can be computed in polynomial time with respect to $d$, when $\p$
is a rational simplex.

Let $L\subseteq V$ be a rational vector subspace.
 Denote the
projected lattice on $V/L$ by $\lattice_L$. Consider the mixed
valuation
$$
S^L(\p,h)= \sum_{y\in \lattice_{V/L}} \int_{\p\cap (y+L)} h(x)dx.
$$
 As shown by Barvinok, and
as we will see here, we can use linear combination of these  mixed
valuations to approximate $S(n\p,h)$ when $n$ is big.

For any polyhedron $\a$, we define the meromorphic function
$S^L(\a,h)(\xi)\in \CM_{\ell}(V^*)$ similarly to $S^L(\a,h)$. For
$\xi\in V^*$ such that the sum converges, we have
$$
S^L(\a,h)(\xi)= \sum_{y\in \lattice_{V/L}} \int_{\a\cap (y+L)}
h(x)e^{\la \xi,x\ra}dm_L(x).
$$
\begin{remark} It is clear that $S^L(\a,h)(\xi)=h(\partial_\xi)\cdot
S^L(\a)(\xi)$.

If  $\p$ is a \emph{polytope}, then $S^L(\p,h)(\xi)$ is regular at
$\xi=0$ and $ S^L(\p,h)(0)= S^L(\p,h)$.
\end{remark}

For a family $\CL$ and a function $L\mapsto\rho(L), L\in\CL,$ we
define
$$
S^{(\CL,\rho)}(\a,h)(\xi)= \sum_{L\in \CL} \rho(L)S^L(\a,h)(\xi).
$$

If  $\p$ is a polytope, and we dilate by $n\in N$, we have again a
quasi-polynomial
$$
S^{(\CL,\rho)}(n\p,h)(0)=\sum_{m=0}^{d+ N} E_m(\CL,\rho,\p,h,n)n
^m.
$$

 We can replace  the quasi-polynomial
$S^{(\CL,\rho)}(n\p,h)(0)$ by   $q$ legal polynomials in the
variable $u$, by splitting  $\N$ into classes modulo $q$. Writing
$n=qu+k$, for  $k=0,\dots,q-1$, we  obtain the polynomial function
of $u$:
$$
S^{(\CL,\rho)}((qu+k) \p,h)(0)= \sum_{m=0}^{d+ N}
E_m^{(k)}(\CL,\rho,\p,h)u ^m.
$$

\bigskip
We briefly recall how the usual  Ehrhart quasi-polynomial of a
polytope can be computed using Brion's theorem. We will then use a
similar method in order to compute efficiently the $r+1$ highest
coefficients only, using Barvinok valuations.

Let $\CV(\p)$ be the set of vertices of $\p$. For each vertex $s$,
let $\c_s$ be the cone of feasible directions of $\p$ at  $s$, so
that the supporting cone at $s$ is $ s+\c_s$. By Brion's theorem
\cite{brion}, we have
$$
S(\p,h)(\xi)= \sum_{s\in \CV(\p)} S(s+\c_s,h)(\xi).
$$
Let $n\in \N $ and consider the dilated polytope $n\p$.
 The supporting  cone at the vertex $ns$ is $ns + \c_s$.
Let $q\in \N$ such that $q \p$ is a lattice polytope and fix $k\in
\N$, $0\leq k\leq q-1$. Let $n=qu+k$. As $qus$ is an integral
point, we have
$$
 S((qu+k)s+ \c_s,h)(\xi)=
 e^{qu\langle
\xi,s\rangle}S^{\CL_s,\rho_s}(ks+\c_s,h)(\xi).
$$
Expanding in powers of $u$, we obtain
$$
S((qu+k)\p,h)(\xi)= \sum_{m\geq 0} u^m \frac{q^m}{m!}\sum_{s\in
\CV(\p)} \langle \xi,s\rangle^m
   \;S(k s+\c_s,h)(\xi).
$$
It follows that for each $m$, the sum of meromorphic functions
$$
\frac{q^m}{m!}\sum_{s\in \CV(\p)} \langle \xi,s\rangle^m
   \;S(k s+\c_s,h)(\xi)
   $$
   is actually analytic. Its value at
$\xi=0$ is obtained by taking the zero degree term. We obtain
$$
S((qu+k)\p,h)(0)=\sum_{m\geq 0} E_m^{(k)}(\p,h)\;u^m,
$$
with
$$
E_m^{(k)}(\p,h)=\frac{q^m}{m!} \sum_{s\in \CV(\p)} \langle
\xi,s\rangle^m
   \;S(k s+\c_s,h)_{[-m]}(\xi).
$$
The right-hand side of this relation, a priori a meromorphic
function of $\xi$, is actually constant. Moreover, we have $S(k
s+\c_s,h)_{[-m]}(\xi)= 0 $ if $m> d+N$, hence the Ehrhart
quasi-polynomial has degree $\leq d+N$.

\bigskip

We apply now Brion's theorem  to $S^{\CL,\rho}(\p,h)(\xi)$. We
obtain
$$
S^{\CL,\rho}(\p,h)(\xi)= \sum_{s\in \CV(\p)}
S^{\CL,\rho}(s+\c_s,h)(\xi).
$$
For reasons to be explained later on, instead of one family $\CL$,
we take a family of subspaces $\CL_s$ for each vertex $s$. Let
$\rho_s: \CL_s\to \Z$ be a function on $\CL_s$. We denote now  by
$(\CL,\rho) $ the map $s\mapsto (\CL_s,\rho_s)$.

We define:
\begin{equation}\label{sommesommets}
\CB^{\CL,\rho}(\p,h)(\xi)= \sum_{s\in
\CV(\p)}S^{\CL_s,\rho_s}(s+\c_s,h)(\xi).
\end{equation}

If the family does not depend on $s$,
$(\CL_s,\rho_s)=(\CL_0,\rho_0)$ {for every vertex $s$}, then, by
Brion's theorem, we have
$$
\CB^{\CL,\rho}(\p,h)(\xi)=S^{\CL_0,\rho_0}(\p,h)(\xi).
$$

We dilate (\ref{sommesommets}). Let $n=qu+k$. We obtain
$$
\CB^{\CL,\rho}((qu+k)\p,h)(\xi)= \sum_{s\in \CV(\p)} e^{qu\langle
\xi,s\rangle}S^{\CL_s,\rho_s}(ks+\c_s,h)(\xi).
$$
Expanding in powers of $u$, we obtain
\begin{equation}\label{sommesommetsdilate}
\CB^{\CL,\rho}((qu+k)\p,h)(\xi)=\sum_{m\geq 0} u^m
E_m^{(k)}(\CL,\rho,\p,h)(\xi)
\end{equation}
with
\begin{equation}\label{ehrhart-of-BCL}
   E_m^{(k)}(\CL,\rho,\p,h)(\xi)=\frac{q^m}{m!}\sum_{s\in
\CV(\p)} \langle \xi,s\rangle^m
   \;S^{\CL_s,\rho_s}(k s+\c_s,h)(\xi).
\end{equation}
If the family does not depend on $s$,
$(\CL_s,\rho_s)=(\CL_0,\rho_0)$ for all vertices, then
$\CB^{\CL,\rho}(\p,h)(\xi) = S^{\CL_0,\rho_0}(\p,h)(\xi)$ is
analytic near $\xi=0$, and so are the coefficients
(\ref{ehrhart-of-BCL}).

On the contrary, if we take a different family $\CL_s$ for each
vertex $s$, the coefficient $ E_m^{(k)}(\CL,\rho,\p,h)(\xi)$  of
$u^m$ in (\ref{sommesommetsdilate}) is no longer analytic near
$\xi=0$, in general. However, the meromorphic function $\xi
\mapsto E_m^{(k)}(\CL,\rho,\p,h)(\xi)$ belongs to
$\CM_{\ell}(V^*)$, thus it has a term of degree $0$  with respect
to $\xi$, given by
\begin{equation}\label{zero-degree-ehrhart-of-BCL}
  E_m^{(k)}(\CL,\rho,\p,h)_{[0]}(\xi)=\frac{q^m}{m!}\sum_{s\in
\CV(\p)} \langle \xi,s\rangle^m  \;S^{\CL_s,\rho_s}(k
s+\c_s,h)_{[-m]}(\xi).
\end{equation}
For a family $(\CL,\rho)$ as described in the next theorem, it turns
out that, for large $m$, this zero-degree part
$E_m^{(k)}(\CL,\rho,\p,h)_{[0]}(\xi)$ is actually analytic, hence
constant, and its value is equal to the $m$-th Ehrhart coefficient
$E_m^{(k)}(\p,h)$  of $S((k+ qu)\p,h)(0)$.

\begin{theorem}\label{betterthanbarvinok}

Let $\p$ be a rational polytope in a rational vector space of
dimension $d$. For each vertex $s$ of the polytope $\p$, let
$\c_s$ be the cone of feasible directions of $\p$ at  $s$, so that
the supporting cone at $s$ is $s+\c_s$. For each vertex $s$, let
$\CL_s$ be a finite family of rational vector subspaces of $V$,
stable under sum, such that $\lin(\f)$ belongs to $\CL_s$ for
every face $\f$ of codimension $r$ of the cone $\c_s$, and let $
\rho_s$ be a patchwork function on $\CL_s$. Let $q\in \N$ such
that $q \p$ is a lattice polytope and fix $k\in \N$, $0\leq k\leq
q-1$. Let $h(x)$ be a homogeneous polynomial of total degree $N$.

Then, for $m\geq d+N-r$, the zero-degree term  $
E_m^{(k)}(\CL,\rho,\p,h)_{[0]}(\xi)$  defined by
(\ref{zero-degree-ehrhart-of-BCL}) is regular near $\xi=0$, hence
constant. Its value  is the coefficient $E_m^{(k)}(\p,h)$ of $u^m$
in the Ehrhart quasi-polynomial
$$
S((k + qu)\p,h)(0)= \sum _{x\in ((k + qu)\p)\, \cap \lattice}h(x)=
\sum_{m=0}^{d+N}u^m E_m^{(k)}(\p,h).
$$
\end{theorem}

\begin{proof}
We first consider the case $h(x)=1$. We have, for every $m\geq 0$,
$$
E_m^{(k)}(\p,1)=\frac{q^m}{m!}\sum_{s\in \CV(\p)} \langle
\xi,s\rangle^m \;S(k s+\c_s)_{[-m]}(\xi)
$$
where the right-hand side is actually a constant function of $\xi$.
For $m
>d-r-1$, we have, by Corollary \ref{lowestdegree},
$$
S(ks+\c_s)(\xi)_{[-m]}= S^{\CL_s,\rho_s}(ks+\c_s)_{[-m]}(\xi).
$$

This proves the theorem when $h(x)=1$. The case  of a non constant
polynomial $h(x)$ is quite similar. If $h(x)= x_1^{N_1}\dots
x_d^{N_d}$, we just have to replace  the meromorphic functions
$S(ks +\c_s)(\xi) $ and $S^{\CL_s,\rho_s}(ks+\c_s)(\xi)$ by their
derivatives under $\partial_{\xi_1}^{N_1}\dots
\partial_{\xi_d}^{N_d}$. $\Box$
\end{proof}

If for each vertex $s$, we take  $\CL_s=\CL$, the full collection
generated by all $r$ codimensional faces of $\p$,
 we obtain Corollary \ref{barvinok} below, that is Barvinok's theorem \cite{barvinokEhrhart},
with an extension to the sum of values of any polynomial $h(x)$ over
the set of integral points of a rational polytope (Barvinok
considers only the case $h(x)=1$).
\begin{corollary}\label{barvinok}
Let $\p\subset V$ be a rational polytope and let $h(x)$ be a
polynomial function on $V$. Let $\CL$ be a finite family of rational
vector subspaces of  $V$, stable under sum. Assume that $\lin(\f)$
belongs to $\CL$ for every  face $\f$ of codimension $r$ of $\p$.
Let $ \rho$ be a patchwork function on  $\CL$ and let
$S^{\CL,\rho}=\sum_{L\in\CL}\rho(L)S^L$. Then the $r+1$ highest
Ehrhart coefficients of $S(t\p, h)(0)$ and $ S^{\CL,\rho}(t\p,
h)(0)$ are equal.
\end{corollary}

The point  in taking a  family $\CL_s$  which depends on the
vertex $s$ lies in the case where   $\p$ is \emph{simplicial}.  In
this case, we can take $\CL_s$ to be just the set of subspaces
$\lin(\f)$, for all faces $\f$ of codimension $\leq r $ of the
supporting cone $\c_s$ at vertex $s$. This family is stable under
sum. Moreover the patchwork function on $\CL_s$ is simple, (Lemma
\ref{rho}) and the computation of the function
$S^{L}(ks+\c_s)(\xi)$, when $L\in \CL_s$, is immediately reduced
(Example \ref{face}) to the computation of a function $S(\a)(\xi)$
for a simplicial cone $\a$ in a rational vector space of dimension
{\bf smaller or equal} than $r$. When $\p$ is a simplex, we obtain
in this way a method for computing the $r+1$ highest Ehrhart
coefficients  for the pair $(\p,h)$.

\section{Local Euler-Maclaurin formula for mixed sums}

Finally in this last section, we  discuss an application of the
existence of the coefficients $\mu^L$ (Theorem \ref{maintheorem})
in the line of \cite{EML}.

Let $\p$ be a rational polytope  in a rational vector space $V$ of
dimension $d$ and let $h(x)$ be a polynomial function on $V$. Let
$L$ be a rational subspace of $V$. Consider the mixed sum
$$
S^L(\p,h)= \sum_{y\in \lattice_{V/L}} \int_{\p\cap (y+L)}
h(x)dm_L(x).
$$

As in \cite{EML}, we associate to the analytic function
$\mu^L(\t(\p,\f))$ a constant coefficients differential operator
(of infinite order) on $V$.

\begin{definition}\label{dpf}
Let $\f$ be a face of $\p$. We denote by $ D^L(\p,\f)$ the
differential operator on  $V$ associated to analytic function
 $\mu^L(\t(\p,\f))$:

 $$D^L(\p,\f)(\partial_\xi)\cdot
 e^{\langle \xi,x\rangle}=\mu^L(\t(\p,\f)))(\xi) e^{\langle
 \xi,x\rangle}.$$
\end{definition}

The operators $D^L(\p,\f)$ are  {\bf local}, that is they depend
only of the transverse cone $\t(\p,\f)$ of $\p$ along $\f$, and
they involve only derivatives in directions orthogonal to the face
$\f$. We can state the following theorem with the same proof as in
\cite{EML}.

\begin{theorem}(Local Euler-Maclaurin formula)\label{LEMLF}
Let $\p$ be a polytope in $V$. For any polynomial function $h(x)$
on $V$, we have
\begin{equation}\label{maclaurin}
S^L(\p,h)=\sum_{\f\in\CF(\p)}\int_\f D^L(\p,\f)\cdot h
\end{equation}
 where the integral on the face $\f$ is taken with respect to the
Lebesgue measure on $<\f>$ defined by the lattice $\lattice\cap
\lin(\f)$.

In particular, for $h=1$, we obtain
\begin{equation}\label{h1}
S^L(\p,1)=\sum_{\f\in\CF(\p)} \mu^L(\t(\p,\f))(0)\vol(\f).
\end{equation}
\end{theorem}

Let us dilate the polytope $\p$ by a non negative integer $n$. If
$\f$ is a face of $\p$,  let $q_{\f}$ be the smallest positive
integer such that $q_{\f}<\f>$ contains integral points. Define
$D(\p,\f,n)=D(n\p,n\f)$, if $n>0$, and
$D(\p,\f,0)=D(q_\f\p,q_\f\f)$. The function $n\mapsto D(\p,\f,n)$
is periodic of period $q_\f$.

\begin{proposition}
Let $\p$ be a rational polytope and $h$ a polynomial function of
degree $N$ on $V$. Then, for any integer $n\geq 0$, we have
\begin{equation}\label{Stph}
S^L(n\p,h)=\sum_{\f\in\CF(\p)}\int_{n\f}D^L(\p,\f,n)\cdot h.
\end{equation}
Furthermore, if $\f\in \CF(\p)$,
  we have
$$
\int_{n\f}D^L(\p,\f,n)\cdot h =\sum_{i=\dim\f}^{\dim\f+
N}E_i(\p,h,\f,n)\,n^i
$$
where the coefficients $E_i(\p,h,\f,n)$ are periodic with period
$q_\f$.

 Hence the Ehrhart coefficients are given by
$$
E_m^L(\p,h,n)= \sum_{\f,{\rm dim \f} \leq m}E_m(\p,h,\f,n).
$$
\end{proposition}

When we apply the last proposition to the function $h(x)=1$, we
obtain
\begin{equation}\label{h1}
S^L(n\p,1)=\sum_{\f\in\CF(\p)} \mu^L(n\t(\p,\f))(0)\vol(\f)
n^{\dim\f}.
\end{equation}
As $\mu^L(n\t(\p,\f))$ is invariant by integral translations, the
function  $\mu^L(n\t(\p,\f))(0)$ is of period $q_\f$.

\thanks{Acknowledgments:
We  thank the  various institutions that helped us to collaborate
  on this work:\ the Research-in-pairs program at the
  Forschungsinstitut Oberwolfach, the University Denis Diderot in Paris and  the Centre
  Laurent Schwartz at Ecole Polytechnique.}

\bigskip

\thanks{{\bf Velleda Baldoni}, Universita  di Roma Tor Vergata,
Dipartimento di Matematica, via della Ricerca Scientifica, 00133
Roma, Italy.}

\thanks{email: baldoni@mat.uniroma2.it}

\thanks{{\bf Nicole Berline},
Ecole Polytechnique,  Centre de math\'ematiques Laurent Schwartz,
91128, Palaiseau, France.}

\thanks{email: berline@math.polytechnique.fr}

\thanks{{\bf Mich{\`e}le Vergne}, Institut de Math\'ematiques de Jussieu, Th{\'e}orie des
Groupes, Case 7012, 2 Place Jussieu, 75251 Paris Cedex 05, France;
 }

\thanks
{Ecole Polytechnique,  Centre de math\'ematiques Laurent Schwartz,
91128, Palaiseau, France.}

\thanks{email: vergne@math.polytechnique.fr}


\begin{thebibliography}{99}
\bibitem{barvinok94}{\bf Barvinok A. I.},\emph{Computing the Ehrhart
polynomial of a convex lattice polytope}, Discrete Comput. geom.
\textbf{12} (1994), 35-48.
\bibitem{barvinokEhrhart}{\bf Barvinok A. I.},{\em Computing the Ehrhart
quasi-polynomial of a rational simplex},  Mathematics of
Computation, {\bf 75} (2006), 1449--1466.
\bibitem{barpom}{\bf Barvinok A. I. and Pommersheim J.},
{\em An algorithmic theory of lattice points in polyhedra}, New
Perspectives in Algebraic Combinatorics (Berkeley,CA, 1996-97),
Math. Sci. Res. Inst. Public {\bf 38}, Cambridge University Press,
Cambridge, (1999), pp 91-147.
\bibitem{EML}{\bf Berline N. and Vergne M.} \emph{Local Euler-Maclaurin formula for polytopes}
(2005), arXiv:math CO/0507256. To appear in Moscow Math. J.
\bibitem{brion}{\bf Brion M.}, {\em Points entiers dans les poly{\`e}dres
convexes}, Ann. Sci. Ecole Norm. Sup. {\bf 21} (1988), 653-663.
\bibitem{latte} \textbf{ De Loera J.A.,  Haws D.,  Hemmecke R.,
 Huggins H.,  Tauzer J.  and
Yoshida R.}, \emph{A User's Guide for LattE v1.1, 2003}, software
package LattE, available at http://www.math.ucdavis.edu/~latte.
\bibitem{Stanley} {\bf Stanley R.} \emph{Enumerative combinatorics}.
 Vol 1 (1997), Cambridge Studies in Advanced Math. \textbf{49}.
\end{thebibliography}
\end{document}